\documentclass[10pt]{amsart}
\pdfoutput=1

\usepackage{cite}

\usepackage{kotex}
\usepackage{lineno,hyperref}
\modulolinenumbers[5]
\usepackage{epstopdf}
\usepackage{pifont}
\usepackage{latexsym}
\usepackage{graphics}
\usepackage{graphicx}
\usepackage{float}
\usepackage[dvips]{xy}
\xyoption{all}
\usepackage{ifpdf}
\usepackage{multirow}
\usepackage{amssymb, amsthm, amsmath}
\usepackage{hhline}
\usepackage{centernot}
\usepackage{verbatim, comment, color}
\usepackage{enumerate}
\usepackage{setspace}
\usepackage{lipsum}% <- For dummy text
\usepackage{subcaption}

\usepackage{mathtools}
\DeclarePairedDelimiter{\ceil}{\lceil}{\rceil}
\DeclarePairedDelimiter{\floor}{\lfloor}{\rfloor}

\newtheorem{theorem}{Theorem} %[section]
\newtheorem{proposition}[theorem]{Proposition}
\newtheorem{lemma}[theorem]{Lemma}
\newtheorem{remark}[theorem]{Remark}

\newtheorem{definition}[theorem]{Definition}

\newcommand{\ba}{\begin{align}}

\newcommand{\be}{\begin{equation}}
\newcommand{\dd}{\rho}
\newcommand{\ee}{\end{equation}}
\newcommand{\bea}{\begin{eqnarray}}
\newcommand{\eea}{\end{eqnarray}}
\newcommand{\barr}{\begin{array}}
\newcommand{\earr}{\end{array}}
\newcommand{\bn}{\begin{enumerate}}
\newcommand{\en}{\end{enumerate}}
\newcommand{\bi}{\begin{itemize}}
\newcommand{\ei}{\end{itemize}}
\newcommand{\bbbm}{\begin{pmatrix}}
\newcommand{\eeem}{\end{pmatrix}}

\newcommand{\bbN}{{\bf N}}

\newcommand{\bbP}{{\bf P}}

\newcommand{\bbR}{{\bf R}}
\newcommand{\bbS}{{\bf S}}

\newcommand{\bbZ}{{\bf Z}}

\newcommand{\cP}{{\cal P}}
\newcommand{\cPf}{\cP_{{\rm fin}}}
\newcommand{\cPu}{\cP^=_{{\rm fin}}}

\newcommand{\cM}{{\cal M}}

\newcommand{\R}{{\mathbf R}}

\newcommand{\al}{\alpha}

\newcommand{\La}{\Lambda}

\newcommand{\ignore}[1]{}{}
\newcommand{\noin}{\noindent}

\newcommand{\nn}{\nonumber}

\newcommand{\p}{{\partial}}

\newcommand{\q}{\quad}

\newcommand{{\QED}}{{\hfill QED} \smallskip}

\newcommand{\spt}{\mathop{\rm spt}}

\renewcommand{\subset}{\subseteq}

\renewcommand{\phi}{\varphi}
\newcommand{\cal}{\mathcal}

\numberwithin{equation}{section}
\numberwithin{theorem}{section}
    
\begin{document}
%\runningtitle{AD}
\title[On a conjecture of Fejes T{\'o}th (1959)]
{On Fejes T{\'o}th's conjectured maximizer for  
 the sum of angles between lines}

\thanks{\em TL is grateful for the support of ShanghaiTech University, and in addition, to the University of Toronto and its Fields Institute for the Mathematical Sciences, where parts of this work were performed.  RM  acknowledges partial support of his research by
Natural Sciences and Engineering Research Council of Canada Grants RGPIN-2015-04383 and 2020-04162.
\copyright 2020 by the authors.
}
\date{\today}

\author{Tongseok Lim and Robert J. McCann}
\address{Tongseok Lim: Krannert School of Management \newline  Purdue University, West Lafayette, Indiana 47907}
\email{lim336@purdue.edu / tlim0213@outlook.com}

\address{Robert J. McCann: Department of Mathematics \newline University of Toronto, Toronto ON Canada}
\email{mccann@math.toronto.edu}

\begin{abstract}
Choose $N$ unoriented lines through the origin of $\R^{d+1}$.  The sum of the angles between these lines is conjectured
to be maximized if the lines are distributed as evenly as possible amongst the coordinate axes of some orthonormal
basis for $\R^{d+1}$.  For $d \ge 2$ we 
embed the conjecture into a one-parameter family of problems,  in which we seek to maximize the sum of
the $\alpha$-th power of the renormalized angles between the lines.  We show the conjecture is equivalent to this same
configuration becoming the {\em unique} optimizer (up to rotations) for all $\alpha>1$.  We establish both the asserted optimality and uniqueness in the limiting case $\alpha =\infty$ of mildest repulsion. The same conclusions extend to
 $N=\infty$, provided we assume only finitely many of the lines are distinct.
\end{abstract}

\maketitle


\section{Introduction}
Choose $N$ unoriented lines through the origin of $\R^{d+1}$.  The sum of the angles between these lines is conjectured
to be maximized if the lines are distributed as evenly as possible amongst the coordinate axes of some orthonormal
basis for $\R^{d+1}$.  When $d=2$ this conjecture dates back to Fejes T\'oth \cite{FT59}.
For $d\ge 2$ it has motivated 
a recent series of works by  Bilyk, Dai, Glazyrin, Matzke, Park, Vlasiuk in different combinations
 \cite{BD19} \cite{BGMV19} \cite{BGMV19-2} \cite{BM19}, and by Fodor, V\'{\i}gh and Zarn\'{o}cz \cite{FVZ16}. 
In this note we 
embed the conjecture into a one-parameter family of problems,  in which we seek to maximize the sum of
the $\alpha$-th power of the renormalized angles between the lines.  Here {\em renormalized} means the angles are rescaled to achieve a maximum value of $1$ when the lines are 
orthogonal, and $\alpha$ parameterizes the effective strength of the repulsion between each
pair of lines.  We show the conjecture is equivalent to this same
configuration becoming the {\em unique} optimizer (up to rotations) for all $\alpha>1$.  We establish both the asserted optimality and uniqueness in the limiting case $\alpha =\infty$ of mildest repulsion.  The optimality extends also to 
$N=\infty$, 
under the assumption that only finitely many of the lines remain distinct.
  
Let $\bbS^d= \{x \in \R^{d+1} \ | \ |x|=1 \}$ be the $d$-dimensional unit sphere, and let $\dd$ be the geodesic distance on $\bbS^d$, i.e. $\dd(x,y) = \arccos (x \cdot y)$. 

Define functions $\La_0 : [0, \pi] \to \R$ and $\La : \bbS^d \times \bbS^d \to \R$ as
\begin{align}
   & \La_0(t)= \frac2\pi \min\{t,\pi -t\} \\ 
   &\La(x,y) = \La_0( \dd(x,y) ),
\end{align}
so that $\La(x,y) \in [0,1]$ is proportional to the (non-obtuse) angle between the lines $x \bbR$ and $y \bbR$ 
in $\bbR^{d+1}$
determined by 
$x$ and $y$ in $\bbS^d$, attaining its maximum value $\La(x,y)=1$ if and only
if $x\cdot y=0$.  Next, we define some measure spaces. Let $N \in \bbN := \{1,2,\ldots\},$
\begin{align*}
\cM(\bbS^d) &= \{ \mu \ | \ \mu \text{ is a nonnegative finite measure on } \bbS^d\}, \\
\cP(\bbS^d) &= \{ \mu  \in \cM(\bbS^d) \ | \ \mu(\bbS^d)=1, \text{ i.e. $\mu$ is a probability}\}, \\
\cP_N(\bbS^d) &= \{ \mu \in \cP(\bbS^d) \ | \  \#[ \spt(\mu) ] \le N \}, \\
\cP_{N}^=(\bbS^d) &= \{ \mu \in \cP_N(\bbS^d) \ | \ \mu = \frac{1}{N}\sum_{i=1}^N \delta_{x_i}, \  x_i \in \bbS^d  \ \forall i  \},  \\
\cPf(\bbS^d) &= \bigcup_{N=1}^\infty \cP_N(\bbS^d), \q {\rm and} \q
\cPu(\bbS^d) =  \bigcup_{N=1}^\infty \cP_{N}^{ =}(\bbS^d).
\end{align*}
Notice that $\cPu(\bbS^d)$ denotes the set of  probability measures with finite support and rational weights on $\bbS^d$.
Then we consider the interaction energy given by $\mu \in \cM(\bbS^d)$,
\be\label{energy1}
E_{1}(\mu) = \frac{1}{2}\iint \La(x,y) d\mu(x)d\mu(y).
\ee

In 1959, Fejes T{\'o}th \cite{FT59} considered the $d=2$ instance of 
the following  problem and several variants
\be\label{FT1}
\text{maximize \ } E_{1}(\mu) \text{ \ over \ } \cP_{N}^{=}(\bbS^d) \ \  \text{for each fixed } N \in \bbN.
\ee
In other words, given $N$, he was interested in the location of $N$ not necessarily distinct points on the sphere which maximizes the sum of their mutual non-obtuse angles.
He conjectured that this energy (discrete sum) is maximized by the periodically repeated copies of the orthonormal basis. More precisely:
\\

\noin{\bf Conjecture 1.} (Fejes T{\'o}th \cite{FT59})\footnote{\cite{BM19} notes the original conjecture was stated for $\bbS^2$ and later stated for all $\bbS^d$ in \cite{P14}.} Given $N \in \bbN$, $\mu$ maximizes \eqref{FT1} if there is an orthonormal basis $v_1,...,v_{d+1}$ of $\R^{d+1}$ and $x_1,...,x_N \in \bbS^d$ such that $\mu = \frac{1}{N}\sum_{i=1}^N \delta_{x_i}$ and $x_i \in \{v_j, -v_j\}$ if $i \equiv j$ mod ${d+1}$.
\\

Here $\delta_x$ is the probability concentrated at $x$, i.e. $\delta_x(A)=1$ if $ x \in A$ and zero otherwise.  We shall also consider a version of the conjecture in which $N$ is unconstrained:
\\

\noin{\bf Conjecture 1'.} A measure $\mu$ maximizes \eqref{energy1} over $\cPf(\bbS^d)$ if there is an orthonormal basis $v_1,...,v_{d+1}$ of $\R^{d+1}$ such that $\mu = \sum_{i=1}^{d+1} (a_i\delta_{v_i} + b_i\delta_{-v_i})$ with $a_i, b_i \ge 0$ and $a_i+b_i = \frac{1}{d+1}$ for all $i=1,...,d+1$.
\\

Essentially, conjectures $1$ (and $1'$) assert that there exists a maximizer of $E_1$ over $\cP_{N}^{=}(\bbS^d)$ (and  over $\cPf(\bbS^d)$ respectively), which is supported on an orthonormal basis. Note that $1'$ implies $1$ if we further restrict $N$ to be divisible by $d+1$.
These conjectures have been settled for the case $d=1$, but for $d \ge 2$ remain open; see e.g. Bilyk and Matzke \cite{BM19}. 
The analogous problem for points (as opposed to lines)
which repel each other throughout the (non-projectivized) sphere
has been analyzed by many authors, as discussed e.g. by Alexander and Stolarsky \cite{AS74}.  
Following ideas developed by 
P\'olya and Szeg\"o \cite{PS31} and others \cite{Bjorck56} \cite{AS74} \cite{BDM18} in that context,
let us consider the kernel $\La^\al$ with $\al>0$, 
the associated quadratic form
\be\label{quadratic}
B_{\al}(\mu, \nu) = \iint \La^\al(x,y) d\mu(x)d\nu(y), \q \mu, \nu \in \cM(\bbS^d),
\ee
and the corresponding energy
\be\label{energy}
E_{\al}(\mu) =  \frac{1}{2} B_\al(\mu,\mu), 
\q \al \in [1, \infty].
\ee
Moreover, let us consider also the limiting case of mildest repulsion: $\La^\infty(x,y)=1$ if $x \cdot y =0$, and zero otherwise. Now we reformulate the conjectures. Notice the following are if and only if statements.
\\

\noin{\bf Conjecture 2.} For each $\al >1$ and $N \in \bbN$, $\mu$ maximizes \eqref{energy} over $\cP_{N}^=(\bbS^d)$ if and only if there is an orthonormal basis $v_1,...,v_{d+1}$ of $\R^{d+1}$ and $x_1,...,x_N \in \bbS^d$ such that $\mu = \frac{1}{N}\sum_{i=1}^N \delta_{x_i}$ and $x_i \in \{v_j, -v_j\}$ if $i \equiv j$ mod ${d+1}$.
\\

\noin{\bf Conjecture 2'.} For each $\al >1$, $\mu$ maximizes \eqref{energy} over $\cPf(\bbS^d)$ if and only if there is an orthonormal basis $v_1,...,v_{d+1}$ of $\R^{d+1}$ such that $\mu = \sum_{i=1}^{d+1} (a_i\delta_{v_i} + b_i\delta_{-v_i})$ with $a_i, b_i \ge 0$ and $a_i+b_i = \frac{1}{d+1}$ for all $i=1,...,d+1$.

\begin{remark} 
\label{rmk1}
Due to the symmetry of $\La_0$ around $t = \pi /2$  we have $\La^\al(x,y) = \La^\al(-x,y)$, which indicates that  the energy \eqref{energy} is invariant under the mass transfer between $x$ and $-x$ on $\bbS^d$. 
Thus the conjecture is more naturally formulated on projective space $\bbR\bbP^d:= \bbS^d/\{\pm\}$,  where the kernel 
$\Lambda^\alpha$ becomes a monotone function of the distance, i.e.\ purely repulsive.
Also, the energy is obviously invariant under the rotation of $\mu$. In view of this, the conjectures 2 and 2' state that for $\al >1$, the maximizer is unique up to rotation and mass exchange between the antipodes. On the other hand, uniqueness is not expected for $\al=1$; e.g. the uniform probability on $\bbS^1$ gives the same $E_1$-energy as $\frac{1}{2}(\delta_x + \delta_y)$ with $x \cdot y =0$.
 \end{remark}

\begin{definition} For $x \in \bbS^d$ we denote $\tilde x = \{x, -x\}$. For $x,y \in \bbS^d$, we denote $x \sim y$ if $\tilde x = \tilde y$. In addition, for $\mu, \nu \in \cM(\bbS^d)$ we write $\mu \equiv \nu$ if  they are equal up to rotation and identification of the antipodes. That is, $\mu \equiv \nu$ if there exists an orthogonal matrix $M \in O(d+1)$ such that $\mu(A \cup -A) = \nu(M(A \cup -A)) $ for every measurable subset $A$ of $\bbS^d$.
 \end{definition}
  
 In this paper we say the optimizer is {\em essentially unique} if $\mu \equiv \nu$ for any optimizers $\mu, \nu$ for a given problem.  In the sequel, we denote Conjecture 1 as {\bf C1}, and {\bf C2}, {\bf C1'}, {\bf C2'} similarly. We begin with a simple observation.
 \begin{proposition}\label{equivalence}
 {\bf C1} and {\bf C2} are equivalent, and so are {\bf C1'} and {\bf C2'}.
  \end{proposition}
\noin{\bf Proof.} Firstly it is clear that {\bf C2} implies {\bf C1} by taking $\al \searrow 1$. For the converse, assume $\bf C1$ holds but {\bf C2} fails for some $\al >1$, so that there is a maximizer $\nu$ of $E_\al$ which is not of the type described in {\bf C1}. Then there exist $x, y \in \spt(\nu)$ such that $x \cdot y \notin \{0, 1, -1\}$. Let $\mu$ be a maximizer of $E_1$ of the type described in {\bf C1}. Then we have the contradiction
 \[
E_1(\mu) =  E_\al(\mu) \le E_\al(\nu) < E_1(\nu) \le E_1(\mu)
 \]
 where the strict inequality follows from $x \cdot y \notin \{0, 1, -1\}$. 
 The equivalence between {\bf C1'} and {\bf C2'} can similarly be understood. \QED
 \\
 
It is clear that if the conclusion of {\bf C2} (resp. {\bf C2'}) holds for some $\al>1$, then the same conclusion also holds  for all $\al'$ satisfying $\al' > \al$. In this paper we resolve {\bf C2} and {\bf C2'} in the mildest case $\al = \infty$, meaning two particles on the sphere interact precisely when they are perpendicular.    
  In a companion work,  by applying a perturbative approach based on Theorem \ref{main} below
    and a bound from \cite{LM20_3},  we are able 
   to generalize the conclusions of that theorem to large finite values of $\al$;
   in that context we are able to extend the maximization to the full set $\cP(\bbS^d)$ \cite{LM20_2}.  
We hope that this may yield an interesting new approach to the Fejes T{\'o}th conjecture in general dimensions.
 
 \begin{theorem}\label{main}
Fix $\alpha=\infty$.
 Then the conclusion of Conjectures 2 and 2' hold in every dimension $d \in \bbN$.
 \end{theorem}

 In section 2 of this paper
the special case $\alpha=\infty$ of {\bf C2} is established
  using an induction on both the dimension $d$ and number of particles $N$
 (see Theorem \ref{main'} for more detailed statement);
   we then deduce from it the more general statement of Theorem \ref{main}  in section~3.

\begin{remark}
After this manuscript was posted on the arXiv,  an elegant alternate proof of 
Theorem 1.4 was provided to us by Dmitriy Bilyk and his collaborators \cite{BGMPV}.
\end{remark}   

   \section{Maximization of $E_\infty$  over $\cP_{N}^=(\bbS^d)$}
   Let $V$ be a $D$-dimensional subspace of $\R^{D+1}$, so that $V \cap \bbS^D$ is a 
$(D-1)$-dimensional unit sphere. Define
\begin{align*}
\cM_{N}^=(\bbS^D) &=N\cP_{N}^=(\bbS^D) = \{ \mu \ | \ \mu = \sum_{i=1}^N \delta_{x_i}, \  x_i \in \bbS^D  \},   \\
\cM_{N,k}^V(\bbS^D) &= \{\mu \in \cM_{N}^=(\bbS^D) \ | \ \mu(V) = k
 \} \ \text{ for } \ k =0, 1,...,N,
 \end{align*}
 \begin{align*}
 E_{D,N} = \max_{\mu \in \cM_{N}^=(\bbS^D)}E_\infty(\mu), \q E_{D,N,k} = \max_{\mu \in \cM^V_{N,k}(\bbS^D)} E_\infty(\mu).
\end{align*}

 \begin{theorem} 
  \label{main'}
  Let $D, N \in \bbN$ and $\al = \infty$. Then:
  \smallskip 
    
 \noin
   (a) The conclusion of {\bf C2} holds. That is, $\mu$ maximizes $E_\infty$ over $\cM_{N}^=(\bbS^D)$ if and only if there is an orthonormal basis $v_1,...,v_{D+1}$ of $\R^{D+1}$ and $x_1,...,x_N \in \bbS^D$ so that $\mu = \sum_{i=1}^N \delta_{x_i}$ and $x_i \sim v_j$ if $i \equiv j$ mod ${D+1}$.
\smallskip 
  
 \noin (b) Assume $k \ge \frac{DN}{D+1}$. Then $\mu = \sum_{i=1}^N \delta_{x_i}$ maximizes $E_\infty$ over   $\cM_{N,k}^V(\bbS^D)$ if and only if $x_i \sim p$ for all $x_i \notin V$ where $ \tilde p = \bbS^D \cap V^\perp$, and there exists an orthonormal basis $v_1,...,v_D$ of $V$ such that either $ \big| \{i \ | \ x_i \sim v_j \} \big| = \ceil[\big]{\frac{k}{D}}$ or $ \big| \{i \ | \ x_i \sim v_j \} \big| = \floor[\big]{\frac{k}{D}}$ for every $j=1,...,D$. 
   \end{theorem}

   \begin{remark}\label{rmk}
   The conclusion of Theorem \ref{main'}(b) holds for $k = \floor{\frac{DN}{D+1}}$ as well, because in this case (along with $k = \ceil{\frac{DN}{D+1}}$) the structure of maximizers over $\cM_{N,k}^V(\bbS^D)$ given in (b) coincides with the  structure of maximizers over $\cM_{N}^=(\bbS^D)$ given in (a), as one can easily check.
      \end{remark}
    
To prove the theorem, we shall first prove the following lemma.
\begin{lemma}\label{comparison}
Let $F_{d,n}=E_\infty(\mu_{d,n})$ and $F_{d,n,k}=E_\infty(\nu_{d,n,k})$ where $\mu_{d,n}, \nu_{d,n,k}$ are the (conjectured) maximizers over $\cM_{n}^{=}(\bbS^d)$ and $\cM_{n,k}^V(\bbS^d)$ respectively as in Theorem \ref{main'}. If $n -1 \ge d \ge 1$ and $\frac{dn}{d+1} \le k \le n-1$ then
\be
\label{inequality}
F_{d, n-k} + F_{d-1,n} < F_{d,n,k}
\ee
 and there is a monotonicity
\begin{align}
\label{mono1}
F_{d,n,k} &> F_{d,n,k+1},\\
\label{mono2}
F_{d,n, \floor{\frac{dn}{d+1}}} &= F_{d,n,\ceil{\frac{dn}{d+1}}}.
\end{align}

\end{lemma}
\noin Of course, $E_{d,n}=F_{d,n}$ and $E_{d,n,k}=F_{d,n,k}$ once the theorem is proved.
\\

  \noin {\bf Proof of Lemma \ref{comparison}.} 
By the structure of maximizers $\mu$ described in Theorem \ref{main'} (a)--(b) respectively, it follows
\begin{align}
\label{eq1}
F_{d, m+1} - F_{d, m} &= m - \floor[\Big]{\frac{m}{d+1}} = \ceil[\Big]{\frac{d}{d+1}m} \ \text{ for every $m$, and} \\
\label{eq2}
F_{d,n,k}  &= F_{d-1,k} +k(n-k).
\end{align} 
Firstly \eqref{mono2} follows from Remark \ref{rmk}. To see \eqref{mono1}, by using \eqref{eq1}, \eqref{eq2} we find \eqref{mono1} is equivalent to $2k-n+1 > \ceil{\frac{d-1}{d}k}$. Notice this is also equivalent to $2k-n \ge \frac{d-1}{d}k$, that is, $k \ge \frac{dn}{d+1}$ which is what we assumed. 

From now on we will prove  \eqref{inequality}. Note that \eqref{inequality} is equivalent to 
\be
\label{ineq}
F_{d, n-k} + F_{d-1,n} - F_{d-1,k}  < k(n-k).
\ee
Firstly, notice $F_{d, n-k}$ is clearly dominated by $E_\infty(\nu)$ where $\nu$ assigns the total mass $n-k$ uniformly onto an orthonormal basis of $\R^{d+1}$, hence 
\be
F_{d, n-k} \le \Big(\frac{(d+1)d}{2}\Big)\Big(\frac{n-k}{d+1}\Big)^2
= \frac{d(n-k)^2}{2(d+1)}. \nn
\ee
Secondly, observe (note that the sum is not void as $n \ge d+1$)
\begin{align*}
F_{d-1,n} - F_{d-1,k} &= \sum_{m=k}^{n-1}  \ceil[\Big]{\frac{d-1}{d}m}  \\
&= e+ \sum_{m=k}^{n-1}{\frac{d-1}{d}m}  \ \text{ where } \  e=e(n,k)  \le n-k \\
&= \frac{(d-1)(n-k)(n+k-1)}{2d} + e.
\end{align*}
It follows
\begin{align*}
&F_{d, n-k} + F_{d-1,n} - F_{d-1,k} \\
& \le  \frac{d(n-k)^2}{2(d+1)} + \frac{(d-1)(n-k)(n+k-1)}{2d} + e \\
&  \le \frac{d(n-k)^2}{2(d+1)} + \frac{(d)(n-k)(n+k-1)}{2(d+1)} + e \\
& = \frac{d(n-k)}{(d+1)}\Big(n-\frac{1}{2}\Big) +e.
\end{align*}
In view of this, \eqref{ineq} will follow if we show
\be \label{ineq1}
e_{{\rm avg}}:= \frac{e}{n-k} < \frac{d}{2(d+1)} + k - \frac{dn}{d+1}.
\ee
Since the ``average error" $e_{{\rm avg}}  \le 1$, \eqref{ineq1} is immediate if $k \ge \ceil{\frac{dn}{d+1}}+1$. So we only need to show \eqref{ineq1} for $k = \ceil{\frac{dn}{d+1}}$, which we henceforth assume. Now let $l := n-k =  \floor{\frac{n}{d+1}}$ so that $(n,k) = \big( (d+1)l+r, dl+r \big)$ for some $r \in \{0,1,...,d\}$. Then since $n= (d+1)l+r$,
\be 
 \frac{d}{2(d+1)} + \ceil[\Big]{\frac{dn}{d+1}} - \frac{dn}{d+1} =  \frac{d}{2(d+1)} + \frac{r}{d+1} = \frac{d + 2r}{2(d+1)}. \nn
\ee
On the other hand, since $k=dl+r$, by definition of $e$, 
\be 
e=\frac{r}{d} + \frac{r+1}{d} + \dots + \frac{d-1}{d} + 0 + \frac{1}{d} + \dots + \frac{d-1}{d} + 0 + \frac{1}{d} + \dots  \nn
\ee
where the sum consists of $l$ terms. By the periodically increasing property of this sequence and the property of arithmetic average, we see that $e_{{\rm avg}} = e/l$ is maximized when $l = d-r$, that is 
\be 
e_{{\rm avg}} = \frac{e}{l}\le \frac{1}{d-r} \Big(\frac{r}{d} + \frac{r+1}{d} + \dots + \frac{d-1}{d}    \Big)
= \frac{d+r-1}{2d}.
\nn
\ee
Now it is easy to check $ \frac{d + 2r}{2(d+1)} > \frac{d+r-1}{2d}$, which proves \eqref{ineq1}. \QED
\\

  \noin {\bf Proof of Theorem \ref{main'}.} Given $D$, the theorem is trivial if $N \le D$. So we assume $N \ge D+1$.  We will proceed by an induction on $N$ and $D$.  Given $d,n \in \bbN$, assume:
  \vspace{2mm}

\noin Induction Hypothesis (IH): Theorem \ref{main'}  (a)--(b) holds for every $N \in \bbN$ if $D \le d-1$, and for every $N \le n-1$ if $D=d$.
  \vspace{2mm}

First we will deduce (a) for $ (D,N)=(d,n)$.  Let $q,r \in \bbZ$ be integers such that $n-1 = q(d+1) +r$ with $r \in \{0,\ldots, d\}$.
By IH(a), we can find a maximizer $\mu_0$ over $\cM_{n-1}^=(\bbS^d)$ such that $\mu_0(V) = dq +r =
\ceil[\big]{\frac{d(n-1)}{d+1}} =  \floor[\big]{\frac{dn}{d+1}}$ and $\mu_0(\{p\}) = q = n-1 -  \floor[\big]{\frac{dn}{d+1}} $. Define $\mu = \mu_0 + \delta_p \in \cM_n^{=}(\bbS^d)$. We will show that $\mu$ is the unique maximizer over $\cM_{n}^{=}(\bbS^d)$.  Observe
\[ E_\infty(\mu) =  E_\infty(\mu_0) + B_\infty(\mu_0, \delta_p)  
 = E_{d, n-1} +   \floor[\Big]{\frac{dn}{d+1}}.
\]
Let $\nu_0$ be any element in $\cM_{n-1}^=(\bbS^d)$ and let $\nu = \nu_0 + \delta_p$. Then
\[ E_\infty(\nu) =  E_\infty(\nu_0) + B_\infty(\nu_0, \delta_p) =  E_\infty(\nu_0) +\nu_0(V).\]
If $E_\infty(\nu_0) = E_{d, n-1}$, then $\nu_0 \equiv \mu_0$ by IH(a) and hence $\nu_0(V) \le 
\floor[\big]{\frac{dn}{d+1}}$. Moreover $\nu_0(V) 
=\floor[\big]{\frac{dn}{d+1}} $ yields $\nu \equiv \mu$. We proved:
 \be
E_\infty(\nu_0) = E_{d, n-1}\text{ and } E_\infty(\nu) \ge E_\infty(\mu) \text{ implies } \nu \equiv \mu.  \nn
\ee
Now suppose $ E_\infty(\nu_0) < E_{d,n-1}$ but still $E_\infty(\nu) \ge E_\infty(\mu)$. This implies $\nu_0(V) >  \floor[\big]{\frac{dn}{d+1}}$. Set $k= \nu_0(V)$ and $k' 
= \floor[\big]{\frac{dn}{d+1}}$, and observe that IH(2) and \eqref{eq2} imply 
\begin{align*}
  E_\infty(\mu) &= E_\infty(\mu_0) + k' = E_{d, n-1, k'} + k' = E_{d-1,k'} + k'(n-k')= F_{d,n,k'}, \\
   E_\infty(\nu) &=  E_\infty(\nu_0) + k \le E_{d, n-1, k} + k = E_{d-1,k} + k(n-k) = F_{d,n,k}.
    \end{align*}
 \eqref{mono1},\eqref{mono2} then gives  $E_\infty(\nu) \le E_\infty(\mu)$, yielding $E_\infty(\nu) = E_\infty(\mu)$ and $k  = \ceil{\frac{dn}{d+1}}$. Now the implied equality $E_\infty(\nu_0) =  E_{d, n-1, k}$ tells us the structure of $\nu_0$ via IH(b), which clearly implies $\mu \equiv \nu$, as in the proof of \eqref{mono2}. This completes the proof of IH(a) for $N=n$.
\\

 Now  we will deduce (b) for $(D,N)=(d,n)$ to complete the induction. Fix $k \ge  \ceil[\big]{\frac{dn}{d+1}}$ and assume $\mu \in \cM_{n,k}^V(\bbS^d)$ achieves the maximum $E_\infty(\mu)=E_{d,n,k}$. We will firstly show that $\mu( \tilde p) > 0$ yields (b). Indeed in this case we can again write $\mu = \mu_0 + \delta_p$, so $E_\infty(\mu)= E_\infty(\mu_0) + k$. Hence $E_\infty(\mu)=E_{d,n,k}$ implies $E_\infty(\mu_0)=E_{d,n-1,k}$. Then IH(b) applied to $\mu_0$ yields $\mu$ must be of the type described in (b), as desired.
  
   Hence from now on we will assume $\mu( \tilde p) = 0$ and show this yields the contradiction $E_\infty(\mu)< E_{d,n,k}$. Let $H= \bbS^d \setminus (V \cup \tilde p)$, and denote $\mu_H, \mu_V$ as the restriction of $\mu$ onto $H$ and $V$ respectively, so that $\mu=\mu_H + \mu_V$. Let $\cal W$ be the set of all $(d-1)$-dimensional subspaces of $V$. Then there exists $W \in \cal W$ such that $ \displaystyle \mu(W) = \max_{X \in \cal W}  \mu(X)$. Now observe
\begin{align*}
 E_\infty(\mu) &= E_\infty(\mu_H) + E_\infty(\mu_V) + B_\infty(\mu_H, \mu_V) \\
 &= E_\infty(\mu_H) + E_\infty(\mu_V) + \sum_{x \in \spt(\mu_H)} \mu(x) B_\infty(\delta_x, \mu_V) \\
 &= E_\infty(\mu_H) + E_\infty(\mu_V) + \sum_{x \in \spt(\mu_H)}\mu(x) \mu(V \cap \{x\}^\perp )  \\
 & \le E_\infty(\mu_H) + E_\infty(\mu_V) + (n-k)\mu(W). 
  \end{align*}
Note that  $E_\infty(\mu_H) \le E_{d,n-k}$. For $E_\infty(\mu_V) + (n-k)  \mu(W)$,  observe
\[   E_\infty(\mu_V) + (n-k)  \mu(W) = E_\infty\big(\mu_V + (n-k) \delta_z\big) \ \text{for some } z \in V \cap W^\perp,
\]
which yields $E_\infty(\mu_V) + (n-k)  \mu(W) \le E_{d-1, n}$ in particular. Then by  IH and Lemma \ref{comparison}, we indeed have 
\begin{align*}
 E_\infty(\mu) &\le E_{d,n-k} + E_{d-1, n} = F_{d,n-k} + F_{d-1, n} < F_{d,n,k} \le E_{d,n,k}.
\end{align*} 
This completes the induction and proves Theorem \ref{main'}.
 \QED
 
    \section{Maximization of $E_\infty$  over $\cPf(\bbS^d)$}
    In this section we show that the same measures maximize the quadratic form $E_\infty(\mu)$
    among all measures $\cPf(\bbS^d)$ having finite support, without imposing any upper bound $N$ on
    the size of the support or rationality on the values of $\mu$.
    We start with a lemma showing that rational quadratic forms on the standard simplex admit rational extrema.
    
\begin{lemma}\label{rationalmax}
Let $A$ be a $(d+1) \times (d+1)$ symmetric matrix with rational entries. Let $\Delta_d$ be the $d$-dimensional standard simplex in $\R^{d+1}$, i.e.
\begin{align*}
\Delta_d &= \{x=(x_1,...,x_{d+1}) \in \R^{d+1} \ | \ \sum_{i=1}^{d+1} x_i =1, \ x_i \ge 0 \  \forall i\}, 
\\ \Delta_d^o &= \{x =(x_1,...,x_{d+1}) \in \Delta_d \ | \ x_i > 0  \text{ for all }i\} 
\end{align*}
its relative interior, and $\p\Delta_d= \Delta_d \setminus \Delta_d^o$ its relative boundary.
Define $f: \Delta_d \to \R$ by $f(x) = \frac{1}{2}x^TAx$.
\smallskip 
  
   \noin
   (a) Then $f$ has maxima and mimima on $\Delta_d$ 
whose entries are rational.\smallskip
  
 \noin
   (b) If $f$ has interior extrema on $\Delta_d$, then it has interior extrema whose entries are rational. That is, suppose there is $x \in \Delta_d^o$ such that $f(x)=M:=\max_{\Delta_d} f$.  Then there is $x' \in \Delta_d^o$ whose entries are rational and $f(x')=M$. The same holds for minima.
   \end{lemma}

\noin{\bf Proof.} We first prove (b), and we will only prove for the maximum and proceed  by an induction on $d$. Let $x \in \Delta_d^o$ attain $f(x)=M$.  Then note that the gradient $\nabla f(x)= Ax$ must be perpendicular to $\Delta_d$ at the optimum $x$, that is, $Ax = ce$ for some $c \in \R$ and $e=(1,1,...,1) \in \R^{d+1}$.

We firstly consider the case $x$ is the {\em only} maximum of $f$ over $\Delta_d$. If $A$ is invertible, then $x=c A^{-1}e$. Now since $A$ has rational entries so does $A^{-1}$. Hence $x$ has rational entries if and only if $c$ is rational. But as $\sum_i x_i =1$, $c$ must be rational and we are done. But if $A$ is not invertible, there is a nonzero $y \in \R^{d+1}$ with $Ay=0$. Now if $y \cdot e =0$, then $x+ty \in \Delta_d$ for  all small $t \in \R$, and we have
\be
g(t) := f(x+ty) = \frac{1}{2}(x+ty)^TA(x+ty) = f(x), \nn
\ee
violating the unique optimality assumed of $x$, thus $y \cdot e \neq 0$. Let $H$ be the hyperplane containing $\Delta_d$. By multiplying a constant we can assume $y \in H$. Now if $y=x$, observe the unique optimality of $x$ implies the system $Az=0$ and $\sum_i z_i =1$ has a unique solution (that is, $x$). Then $x$ will be found via elementary row operations on this system, and the rationality of $A$ implies rationality of $x$, yielding the lemma.  On the other hand if $y \neq x$, let $v= y-x$ and $z = x +tv = (1-t)x+ty$. Then $z \in \Delta_d$ for  all small $t \in \R$, and $g(t)=f(z)=(1-t)^2f(x)$, again violating the unique optimality of $x$.

It remains to consider the case $x \in \Delta_d^o $ is not the unique maximum, so there is $x' \in \Delta_d$, $x' \neq x$ and $f(x')=M$. Consider the line $L=\{z \ | \ z=(1-t)x +tx',  \ t \in \R\}$, and the two distinct points $\{y, y'\} = L \cap \p\Delta_d$. Observe then there exist two unique subsimplexes $\Delta_k, \Delta_{k'}$ of dimensions $k$ and $k'$ respectively, such that $y, y'$ lie in the relative interiors of $\Delta_k, \Delta_{k'}$, denoted by $\Delta_k^o, \Delta_{k'}^o$. Here a single point is its own relative interior, by convention.
Notice $\Delta_k^o \cap \Delta_{k'}^o = \emptyset$. Now observe that $f$ is constant on $L$, since with $v=x'-x$ we have 
\be
g(t) := f(x+tv) = \frac{1}{2}(x+tv)^TA(x+tv) = f(x) +t^2f(v)  \nn
\ee
where we used the fact that $v \cdot Ax =v \cdot ce=0$. Moreover $g(0)=g(1)$ implies $f(v)=0$, thus $g$ is constant, hence $f(y)=f(y')=M$. Then since $k, k' < d$, the induction hypothesis is applicable and there exist $z \in \Delta_k^o$ and $z' \in \Delta_{k'}^o$ such that $f(z)=f(z')=M$ and the entries of $z, z'$ are rational. Now since the quadratic form $f$ is constant on the lines $\overline{xz}$ and $\overline{xz'}$ as shown above, $f=M$ on the line $\overline{zz'}$. Finally, notice $\overline{zz'} \cap \Delta_d^o \neq \emptyset$; otherwise the segment $\overline{zz'}$ is in $\p \Delta_d$, thus in a subsimplex $\Delta_{d-1}$ of dimension $d-1$. Then $\{z,z' \} \subset  \Delta_{d-1} $, yielding $\{y,y' \} \subset \Delta_{d-1} $ and hence $\overline{yy'} \subset \Delta_{d-1}$. But this contradicts to the fact $x \in \overline{yy'}$. Hence $\overline{zz'} \cap \Delta_d^o \neq \emptyset$, yielding $\{(1-t)z+ tz' \ | \ t \in (0,1) \} \subset \Delta_d^o$. We conclude any rational $t \in (0,1)$ yields the rational maximum $(1-t)z+ tz' \in\Delta_d^o$.

Now part (a) is immediate from part (b) since any (extremal) point $x \in \Delta_d$ belongs to the relative interior of a subsimplex $\Delta_k^o$.
\QED

\begin{theorem}\label{main2} A measure
$\mu$ maximizes $E_\infty$ over $\cPf(\bbS^d)$ if and only if there is an orthonormal basis $v_1,...,v_{d+1}$ of $\R^{d+1}$ such that $\mu \equiv \frac{1}{d+1}\sum_{i=1}^{d+1} \delta_{v_i}$.
\end{theorem}

 \noin{\bf Proof.} Let $\mu^*$ be a conjectured (essentially unique) maximizer, so that $E_\infty(\mu^*)=\frac{(d+1)d}{2}\frac{1}{(d+1)^2} = \frac{d}{2(d+1)}=:M^*$. Choose any $\mu \in \cPf(\bbS^d)$. We will show $E_\infty(\mu) \le E_\infty(\mu^*)$, and the inequality is strict unless $\mu \equiv \mu^*$.
 
 To this end, set $A= \spt(\mu)$ and $\cP(A)=\{ \mu \in \cPf(\bbS^d) \ | \ \spt(\mu) \subset A\}$ and
$M_A:=\max_{\cP(A)} E_\infty$. 
 We separate two possible cases, according to whether or not equality holds in $E_\infty(\mu) \le M_A$.

In the first case where $E_\infty(\mu)=M_A$, notice that Lemma \ref{rationalmax}(b) implies there exists $\nu \in \cP(A)$ such that $\spt(\nu)=A$, $E_\infty(\nu) = M_A$ and $\nu(x)$ is rational for every $x \in A$. Then Theorem \ref{main'}(a) implies $M_A \le M^*$ thanks to the rationality, and furthermore equality holds if and only if $\nu \equiv \mu^*$ (here we use the standard fact that the interaction energy $\sum_{i \ne j} x_ix_j$ is uniquely maximized over the hyperplane $\{ x \in \R^d \ | \ \sum_i x_i=1\}$ by the uniform probability vector;  this can, for example, be seen as a consequence
of the Perron-Frobenius theorem). Hence in particular   $M_A=M^*$ yields $\spt(\mu)  = \spt(\nu)\subset \bigcup_{i=1}^{d+1} \{u_i, -u_i\}$ for some orthonormal basis $u_1,...,u_{d+1}$, and this clearly implies $\mu \equiv \mu^*$ again by the same standard fact.

In the second case $E_\infty(\mu)<M_A$, choose any  $\nu \in \cP(A)$ attaining $E_\infty(\nu) = M_A$. Let $B= \spt(\nu)$ and define $\cP(B)$ and $M_B$ as above. Then by the above argument we have $M_B \le M^*$.  But notice $M_B=E_\infty(\nu) =M_A$, which yields $E_\infty(\mu) < M_A = M_B \le M^*$. 

We conclude $\mu^*$ is the essentially unique maximizer on $\cPf(\bbS^d)$.
 \QED

\noin {\bf Acknowledgement.} This is a post-peer-review, pre-copyedit version of an article published in Applied Mathematics \& Optimization. The final authenticated
version is available online at: \url{https://doi.org/10.1007/s00245-020-09745-5}.

\end{document}